\documentclass[12pt]{amsart}

\usepackage{amssymb,amsthm,amsmath}
\usepackage[utf8]{inputenc}
\usepackage{enumerate}

\newcommand{\R}{\mathbb{R}}
\newcommand{\1}{\textbf{1}}
\newcommand{\dd}{\mathrm{d}}
\newcommand{\E}{\mathbb{E}}

\newcommand{\fun}[3]{#1\colon #2 \longrightarrow #3}

\DeclareMathOperator{\sgn}{sgn}

\newtheorem*{thm*}{Theorem}

\theoremstyle{definition}
\newtheorem{remark}{Remark}

\title{A note on the tensor product of two random unitary matrices}

\author{Tomasz Tkocz}

\address{Mathematics Institute, University of Warwick, Coventry CV4 7AL, UK}
\email{t.tkocz@warwick.ac.uk}

\begin{document}

\begin{abstract}
In this note we consider the point process of eigenvalues of the tensor product of two independent random unitary matrices of size $m \times m$ and $n \times n$. When $n$ becomes large, the process behaves like the superposition of $m$ independent sine processes. When $m$ and $n$ go to infinity, we obtain the Poisson point process in the limit.
\end{abstract}

\maketitle

\section{Introduction}

In quantum mechanics the time evolution of two noninteracting subsystems can be described by an operator $e^{itH} \otimes e^{itH'}$, where $H$ and $H'$ are Hamiltonians of the subsystems (see e.g. chapters 2.2 and 3.1 in \cite{BP}). In applications, the unitary operator $e^{itH}$, which is \emph{a priori} complicated, is replaced by a random unitary matrix, to make a model tractable. This powerful idea goes back to E. Wigner. Here by a $n \times n$ random unitary matrix we mean a matrix drawn according to the Haar measure on the unitary group $U(n)$. From this point of view it seems natural to study asymptotic local properties of spectra of the tensor product $A_m \otimes B_n$ of two independent $m \times m$ and $n \times n$ random unitary matrices, to which this short note is devoted. The note, in a sense, continues the investigations commenced in \cite{T}.

Some preliminaries are presented in the rest of this section, and the main result is stated. The proofs are provided in the next section. The last section is devoted to some concluding remarks concerning the tensor product of more than two matrices.

\subsection{Background and notation}

For a simple point process $\tau$ on $\R$ we denote its \emph{$k$-th correlation function}, when it exists, by $\rho_\tau^{(k)}$ (for the definitions see e.g. \cite{HKPV}). Let us introduce three point processes $\Pi$, $\Sigma$, and $\Xi_n$. By $\Pi$ we shall denote \emph{the Poisson point process} on $\R$ for which $\rho_{\Pi}^{(k)} \equiv 1$ for all $k$. By $\Sigma$ we shall denote \emph{the sine point process} on $\R$ which has the correlation functions
\begin{equation}
\label{eq:defsinecorr}
\rho_{\Sigma}^{(k)}(x_1, \ldots, x_k) = \det\left[ Q(x_i,x_j) \right]_{i,j = 1}^k,
\end{equation} 
where \emph{the sine kernel} $Q(x,y) = q(x-y)$ and $q$ reads as follows
\begin{equation}
\label{eq:defsinekernel}
q(u) = \frac{\sin(\pi u)}{\pi u}.
\end{equation}
Given a $n \times n$ random unitary matrix with eigenvalues $e^{i\xi_1}, \ldots, e^{i\xi_n}$, where $\xi_i \in [0,2\pi)$ are \emph{eigenphases}, we define the point process $\Xi_n = \{\xi_1, \ldots, \xi_n\}$. It is well known that this process is determinantal with the kernel $S_n(x,y) = s_n(x-y)$, where
\begin{equation}
\label{eq:defkerS_n}
s_n(u) = \frac{1}{2\pi}\frac{\sin\left(\frac{nu}{2}\right)}{\sin\left( \frac{u}{2} \right)},
\end{equation}
i.e.,
\begin{equation}
\label{eq:defCUEcorr}
\rho_{\Xi_n}^{(k)}(x_1, \ldots, x_k) = \det\left[ S_n(x_i, x_j) \right]_{i,j = 1}^k.
\end{equation}
Since $\frac{2\pi}{n}s_n\left( \frac{2\pi}{n}u \right) \xrightarrow[n \to \infty]{} q(u)$, when $n$ becomes large, the process $\frac{n}{2\pi}(\Xi_n-\pi)$ of the rescaled eigenphases of the $n \times n$ random unitary matrix locally behaves as the sine process $\Sigma$.

By \emph{superposition} of two simple point processes $\Psi = \{\psi_1, \ldots, \psi_M\}$, $\Phi = \{\phi_1, \ldots, \phi_N\}$, $M, N \leq \infty$ we mean the union $\Psi \cup \Phi = \{\psi_1, \ldots, \psi_M, \phi_1, \ldots, \phi_N\}$.

\subsection{Results}

Given two independent $m \times m$ and $n \times n$ random unitary matrices $A$ and $A'$ we get two independent point processes of their eigenphases $\Xi_m = \{\xi_1, \ldots, \xi_m\}$ and $\Xi_n' = \{\xi'_1, \ldots, \xi'_n\}$ respectively. We define the point process $\Xi_m \otimes \Xi_n'$ of the eigenphases of the matrix $A \otimes A'$ as
\[ \Xi_m \otimes \Xi_n' = \{\xi_i + \xi'_j \ \textrm{mod} 2\pi, \ i = 1,\ldots, m, j = 1, \ldots, n \}.\]
It has been recently shown \cite[Theorem 1]{T} that the process $\frac{n^2}{2\pi}(\Xi_n \otimes \Xi_n')$ behaves locally as the Poisson point process on $\R_+$. We refine this result and investigate what happens when $n$ becomes large with $m$ being fixed, or when both $m$ and $n$ becomes large but not necessarily $m = n$.
\begin{thm*}
\label{thm:main}
Let $\Xi_m$ and $\Xi_n'$ be point processes of eigenphases of two independent $m \times m$ and $n \times n$ random unitary matrices. Let $\Sigma_1, \ldots, \Sigma_m$ be independent sine processes and let $\Pi$ be a Poisson process on $\R$. Then for each $k \leq n$ the $k$-th correlation function of the process $\Xi_m \otimes \Xi_n'$ exists and
\begin{enumerate}[(a)]
\item\label{thm(a)}
$\rho_{\frac{mn}{2\pi}(\Xi_m \otimes \Xi_n' - \pi)}^{(k)} \xrightarrow[n \to \infty]{} \rho_{m\Sigma_1 \cup \ldots \cup m\Sigma_m}^{(k)}$,
\item\label{thm(b)}
$\rho_{\frac{mn}{2\pi}(\Xi_m \otimes \Xi_n' - \pi)}^{(k)} \xrightarrow[m,n \to \infty]{} \rho_{\Pi}^{(k)}$,
\end{enumerate}
uniformly on all compact sets in $\R^k$.
\end{thm*}

\begin{remark}[Weak convergence]\label{rem:weakconv}
According to \cite{HKPV}, by a point process on $\R$ we mean a random variable with values in the metric space $\mathcal{M}(\R)$ of $\sigma$-finite Borel measures on $\R$ (counting measures correspond to locally finite subsets of $\R$) endowed with the topology generated by the functions $\mu \mapsto \int f\dd \mu$ for continuous, compactly supported $f$. We say that a sequence of point processes $(\tau_n)$ converges \emph{in distribution} to a point process $\tau$ if the law $\nu_n$ of $\tau_n$ converges weakly to that of $\tau$, say $\nu$, in the space $\mathcal{M}_1(\mathcal{M}(\R))$ of probability measures on $\mathcal{M}(R)$, i.e. $\int f \dd \nu_n \to \int f \dd \nu$ for any bounded continuous function on $\mathcal{M}(\R)$. Clearly, these integrals can be expressed using correlation functions, hence the theorem implies the convergence in distribution of the considered point processes.
\end{remark}

\begin{remark}[Heuristic behind (\ref{thm(a)})]\label{rem:heuristic}
In view of the mentioned theorem from \cite{T} result (\ref{thm(b)}) should not be surprising. Neither is (\ref{thm(a)}) as in the simplest case $m = 2$ we have
\begin{align*} 
\Xi_2 \otimes \Xi_n' = &\{ \xi_1 + \xi'_1 \ \textrm{mod} 2\pi, \ldots, \xi_1 + \xi'_n \ \textrm{mod} 2\pi \} \\
&\cup \{ \xi_2 + \xi'_1 \ \textrm{mod} 2\pi, \ldots, \xi_2 + \xi'_n \ \textrm{mod} 2\pi \}.
\end{align*}
After shifting and rescaling we end up with two families of the rescaled eigenphases of a $n \times n$ random unitary matrix which differ roughly by a large shift $\frac{n}{2\pi}(\xi_1 - \xi_2)$ which is independent on the matrix. That makes the families independent and in the limit, according to $\rho_{\frac{n}{2\pi}(\Xi_n - \pi)}^{(k)}\xrightarrow[n\to\infty]{}\rho_{\Sigma}^{(k)}$, they look like sine processes. \qedsymbol
\end{remark}

\begin{remark}[Superposition of many sine processes becomes a Poisson point process]\label{rem:superposition}
Notice that for any independent copies $\Phi_1, \ldots, \Phi_m$ of a point process $\Phi$ we have
\[ \rho_{\Phi_1\cup\ldots\cup\Phi_m}^{(k)}(x_1, \ldots, x_k) = \sum_{p=1}^{m \wedge k} \sum_{\pi \in \mathfrak{S}(k,p)} \frac{m!}{(m-p)!} \prod_{j=1}^p \rho_\Phi^{(\sharp \pi_j)} ((x_i)_{i \in \pi_j}),\]
where $\mathfrak{S}(k, p)$ is the collection of all partitions into $p$ nonempty pairwise disjoint subsets of the set $\{1, \ldots, k\}$. By this we mean that if $\pi$ is such a partition then $\pi = \{\pi_1, \ldots, \pi_p\}$, where $\pi_q = \{\pi(q,1),\ldots, \pi(q,\sharp \pi_q)\}$ is the $q$-th block of the partition $\pi$.

Along with the fact that if we rescale, $\rho_{\lambda \Phi}^{(k)}(x)$ becomes $\frac{1}{\lambda^k}\rho_\Phi^{(k)}\left(\frac{1}{\lambda}x\right)$, the previous observation yields
\begin{equation}
\label{eq:remsuperposofsines}
\rho_{m\Sigma_1\cup\ldots\cup m \Sigma_m}^{(k)}(x) = \sum_{p=1}^{m \wedge k} \sum_{\pi \in \mathfrak{S}(k,p)} \frac{1}{m^k}\frac{m!}{(m-p)!} \prod_{j=1}^p \rho_\Sigma^{(\sharp \pi_j)} \left(\frac{1}{m}(x_i)_{i \in \pi_j}\right).
\end{equation}
When $m$ goes to infinity we thus get
\[ \lim_{m \to \infty} \rho_{m\Sigma_1\cup\ldots\cup m \Sigma_m}^{(k)}(x) = \lim_{m \to \infty} \prod_{j=1}^p \rho_\Sigma^{(1)} \left(\frac{1}{m}(x_i)_{i \in \pi_j}\right) = 1 = \rho_\Pi^{(k)}.\]
It retrieves the special case of a quite expected phenomenon put forward in \cite{CD}. Namely, the authors say ``[...] a 
Poisson process can be viewed as an infinite superposition of determinantal
or permanental point processes'' (see Theorem 4 therein and the two preceding paragraphs).
Regarding Theorem (\ref{thm(a)}) that implies
\[ \lim_{m\to\infty}\lim_{n\to\infty}\rho_{\frac{mn}{2\pi}(\Xi_m \otimes \Xi_n' - \pi)}^{(k)} = 1.\]
Note that in the second part of the theorem we establish a stronger statement, that letting the dimensions of two independent random unitary matrices to infinity reduces all the correlations in their tensor product.\qedsymbol
\end{remark}

\section{Proofs}

For the sake of convenience, let us recall a few basic facts which shall be frequently used. 

Note the following easy estimate (for the definition see \eqref{eq:defkerS_n})
\begin{equation}
\label{eq:easyboundS_n}
\sup_{x \in \R}\left|\frac{2\pi}{n}s_n(x)\right| = 1.
\end{equation}
Combined with Hadamard's inequality (see e.g. (3.4.6) in \cite{AGZ}), it allows us to bound the correlation functions,
\begin{equation}
\label{eq:corrbound}
\sup_{x \in \R^k} \rho^{(k)}_{\Xi_n}(x) \leq k^{k/2}\| s_n \|_\infty^k = \frac{k^{k/2}}{(2\pi)^k}n^k.
\end{equation}

\subsection{Proof of Theorem (\ref{thm(a)})}

Let $\Theta_{m,n} = \frac{mn}{2\pi}(\Xi_m \otimes \Xi_n' - \pi)$. Fix a natural number $k$. Since we will let $n$ go to infinity, we may assume that $k \leq n$. First we show that there exists functions $\fun{\rho^{(k)}_{\Theta_{m,n}}}{\R^k}{[0,\infty)}$ so that for any bounded and measurable function $\fun{f}{\R^k}{\R}$ we have
\[ \E \sum f(\theta_1, \ldots, \theta_k) = \int_{\R^k} f(x)\rho^{(k)}_{\Theta_{m,n}}(x) \dd x, \]
where the summation is over all ordered $k$-tuples $(\theta_1, \ldots, \theta_k)$ of distinct points of $\Theta_{m,n}$. This will prove that $\rho^{(k)}_{\Theta_{m,n}}$ are the correlation functions of $\Theta_{m,n}$. Then we will deal with the limit when $n \to \infty$.

Fix $f$. Since for each $s = 1, \ldots, k$, $\theta_s = \frac{mn}{2\pi}(\xi_{i_s} + \xi'_{j_s} \ \textrm{mod} 2\pi - \pi)$ for some $i_s \in \{1,\ldots, m\}$, $j_s \in \{1,\ldots,n\}$ we can write
\[ \E \sum f(\theta_1, \ldots, \theta_k) = \E\sum_{\substack{i \in \{1,\ldots,m\}^k \\ j \in \{1,\ldots,n\}^k}} f\left( \left( \frac{mn}{2\pi}(\xi_{i_s} + \xi'_{j_s} \ \textrm{mod} 2\pi - \pi) \right)_{s=1}^k \right) ,\]
where the second sum is subject to $k$-tuples $i$, $j$ such that the pairs $(i_1, j_1), \ldots, (i_k,j_k)$ are pairwise distinct. For sure it happens when all the $j_s$'s are distinct. Call these choices of $i$ and $j$ \emph{good} and the rest \emph{bad}. So
\[ \E \sum_{i,j} f = \E \sum_{\textrm{good } i,j} f + \E\sum_{\textrm{bad } i,j}f.\] 
First we handle the \emph{good} sum. Some $i_s$'s may overlap and we will control it using partitions of the set $\{1, \ldots, k\}$ into $p \leq k \wedge m$ nonempty pairwise disjoint subsets (see Remark \ref{rem:superposition} for the notation) so that $i_s = i_t$ whenever $s$ and $t$ belong to the same block of a partition. We have
\[ \E\sum_{\textrm{good } i,j} f = \sum_{p=1}^{k \wedge m} \sum_{\pi \in \mathfrak{S}(k,p)} \E \sum_{\substack{\textrm{distinct} \\ i_{\pi(1,1)}, \ldots, i_{\pi(p,1)}}} \sum_{\substack{ \textrm{distinct} \\ j_1, \ldots, j_k}} f.\]
The sums over $i$'s and $j$'s have been separated. Therefore taking advantage of independence as well as recalling definitions of the $p$-th and $k$-th correlation functions of $\Xi_m$ and $\Xi_n'$ we find
\begin{align*}
\E\sum_{\textrm{good } i,j} f = \sum_{p, \pi} \int_{[0,2\pi]^p}\int_{[0,2\pi]^k}& f\left( \left( \frac{mn}{2\pi}(x_{\pi(s)} + y_s \ \textrm{mod} 2\pi - \pi) \right)_{s=1}^k \right) \\ &\rho^{(p)}_{\Xi_m}(x_1, \ldots, x_p)\rho^{(k)}_{\Xi_n'}(y_1,\ldots, y_k) \dd x_1 \ldots \dd x_p \dd y_1 \ldots \dd y_k,  
\end{align*}
where we note $\pi(s) = q  \ \Longleftrightarrow  \ s \in \pi_q$. Finally, we need to address the technicality concerning the addition $\textrm{mod} 2\pi$. Keeping in mind that we integrate over $[0,2\pi]^p$ and $[0,2\pi]^k$ we consider for $\eta \in \{0,1\}^k$ the set
\begin{align*}
U_\eta = \bigg\{ x \in [0,2\pi]^p, y \in [0,2\pi]^k; \ \forall s \leq k \ &x_{\pi(s)} + y_s < 2\pi \textrm{ if } \eta_s = 0, \textrm{ and } \\ &x_{\pi(s)} + y_s \geq 2\pi \textrm { if } \eta_s = 1 \bigg\}.
\end{align*} 
Then on $U_\eta$ we have $x_{\pi(s)} + y_s \textrm{ mod} 2\pi = x_{\pi(s)} + y_s - 2\pi\eta_s$, thus changing the variables on $U_\eta$ so that $z_s  = \frac{mn}{2\pi}(x_{\pi(s)} + y_s - 2\pi\eta_s - \pi)$ we get
\begin{align*}
\E\sum_{\textrm{good } i,j} f = \int_{\R^k} f(z) \Bigg( \sum_{p,\pi,\eta} \1_{W_\eta}(z)\int_{[0,2\pi]^p}\1_{V_\eta}(x)\rho^{(p)}_{\Xi_m}(x)\left( \frac{2\pi}{mn} \right)^k\rho^{(k)}_{\Xi_n'}(y(z,x)) \dd x \Bigg) \dd z,
\end{align*}
where $y_s(z,x) = \frac{2\pi}{mn}z_s - x_{\pi(s)} + 2\pi\eta_s + \pi$,
\[ V_\eta = \left\{ x \in \R^p; \ \forall s \leq k \ \frac{2\pi}{mn}z_s + 2\pi\eta_s - \pi \leq x_{\pi(s)} \leq \frac{2\pi}{mn}z_s + 2\pi\eta_s + \pi \right\},\]
and
\[ W_\eta = \left\{ z \in \R^k; \ \forall s \leq k \ z_s \leq mn/2 \textrm{ if } \eta_s = 0, \textrm{ and } z_s \geq -mn/2 \textrm { if } \eta_s = 1 \right\}.\]
Summarizing, we have just seen that the correlation function $\rho^{(k)}_{\Theta_{m,n}}(z)$ takes on the form
\begin{equation}\label{eq:rho}
\rho^{(k)}_{\Theta_{m,n}}(z) = \sum_{p,\pi,\eta} \1_{W_\eta}(z)\int_{[0,2\pi]^p}\1_{V_\eta}(x)\rho^{(p)}_{\Xi_m}(x)\left( \frac{2\pi}{mn} \right)^k\rho^{(k)}_{\Xi_n'}(y(z,x)) \dd x + B_{m,n}(z),
\end{equation}
where the term $B_{m,n}$ corresponds to the sum over bad indices $\E \sum_{\textrm{bad } i, j} f$. By the same kind of reasoning we show that roughly
\begin{align*} 
B_{m,n}(z) = \sum_{p=1}^k\sum_{q=1}^{k-1}\sum_{\substack{\pi \in \mathfrak{S}(k,p) \\ \tau \in \mathfrak{S}(k,q)}} \sum_\eta \1_{\tilde W_\eta}(z) \left(\frac{2\pi}{mn}\right)^k \int_{[0,2\pi]^{p+q-k}} &\1_{\tilde V_\eta}(x) \rho^{(p)}_{\Xi_m}(\tilde x(z,x)) \\ &\rho^{(q)}_{\Xi_n'}(\tilde y(z,x)) \dd x,
\end{align*}
where the sums are over appropriate partitions and $\tilde W_\eta$, $\tilde V_\eta$ are suitable sets which appear after changing the variables. Now, by \eqref{eq:corrbound},
\begin{equation}
\label{eq:prodbound}
 \|\rho^{(p)}_{\Xi_m}\cdot \rho^{(q)}_{\Xi_n'}\|_\infty \leq \frac{p^{p/2}q^{q/2}}{(2\pi)^{p+q}}m^pn^q, 
\end{equation}  
so
\[ B_{m,n}(z) \leq C_k\frac{1}{n},\]
where the constant $C_k$ depends only on $k$ (roughly, it equals the number of summands times $k^k$). Hence, when taking $n \to \infty$ we will not have to take care about $B_{m,n}$.

Let us have a look at \eqref{eq:rho} and compute now the limit of the first term when $n \to \infty$. We observe that $\1_{W_\eta} \to 1$ pointwise on $\R^k$. Moreover, $\sum_\eta \1_{V_\eta} \to \1_{[0,2\pi)^p}$, and $\1_{V_\eta} \to 0$ for $\eta$ such that $\eta_s \neq \eta_t$ but $\pi(s) = \pi(t)$ for some $s \neq t$. Thus we consider only $\eta$'s such that $\eta_s = \eta_t$ whenever $\pi(s) = \pi(t)$ and then the following simple observation
\begin{equation}
\label{eq:simpleobserv}
\frac{2\pi}{mn}s_n\left( \frac{2\pi}{mn}u + v \right) \xrightarrow[n\to\infty]{} \begin{cases} 0, & v \neq 0 \\ \frac{1}{m}q\left(\frac{u}{m}\right), & v = 0  \end{cases}
\end{equation}
yields for all these $\eta$'s,
\begin{align*}
\left( \frac{2\pi}{mn} \right)^k\rho_{\Xi_n'}^{(k)}(y) &= \det \left[ \frac{2\pi}{mn}s_n\left( \frac{2\pi}{mn}(z_s - z_t) + 2\pi(\eta_s - \eta_t) + x_{\pi(t)} - x_{\pi(s)} \right) \right]_{s,t=1}^k \\
&\xrightarrow[n\to\infty]{} \prod_{j=1}^p \det\left[ \frac{1}{m}q\left( \frac{z_s - z_t}{m} \right) \right]_{s,t \in \pi_j} = \frac{1}{m^k}\prod_{j=1}^p \rho_\Sigma^{(\sharp \pi_j)}\left( \frac{1}{m}(z_i)_{i \in \pi_j} \right).
\end{align*}
By estimate \eqref{eq:corrbound}, $\left( \frac{2\pi}{mn} \right)^k\rho_{\Xi_n'}^{(k)}(y)$ is bounded by $k^{k/2}/m^k$, so the integrand in \eqref{eq:rho} can be simply bounded. Thus by Lebesgue's dominated convergence theorem
\begin{align*}
 \rho_{\Theta_{m,n}}^{(k)}(z) \xrightarrow[n\to\infty]{} \sum_{p,\pi} \frac{1}{m^k}\prod_{j=1}^p \rho_\Sigma^{(\sharp \pi_j)}\left( \frac{1}{m}(z_i)_{i \in \pi_j} \right) \cdot \int_{[0,2\pi]^p} \rho_{\Xi_m}^{(p)}(x) \dd x.
\end{align*}
For any $p \leq m$ the integral $\int_{[0,2\pi)^p} \rho_{\Xi_m}^{(p)}(x) \dd x$ just equals  $m!/(m-p)!$. Consequently, we finally obtain
\[ \rho_{\Theta_{m,n}}^{(k)}(z_1, \ldots, z_k) \xrightarrow[n\to\infty]{} \sum_{p,\pi} \frac{1}{m^k}\frac{m!}{(m-p)!}\prod_{j=1}^p\rho_\Sigma^{(\sharp \pi_j)}\left( \frac{1}{m}(z_i)_{i \in \pi_j} \right).\]
In view of \eqref{eq:remsuperposofsines} this completes the proof. \qedsymbol

\subsection{Proof of Theorem (\ref{thm(b)})}

Fix a point $z = (z_1, \ldots, z_k) \in \R^k$. We let $m$ and $n$ tend to infinity and want to prove that $\rho_{\Theta_{mn}}^{(k)}(z)$ tends to 1. Recall \eqref{eq:rho} and notice that due to estimate \eqref{eq:prodbound} all the terms with $p \leq k-1$ are bounded above by $C_k/m$, so we can write
\begin{align*}
\rho_{\Theta_{m,n}}^{(k)}(z) = O\left(\frac{1}{m}+\frac{1}{n}\right) +   \sum_\eta \1_{W_\eta}(z) \int_{[0,2\pi]^k}\1_{V_\eta}(x) \left( \frac{2\pi}{mn} \right)^k  \rho_{\Xi_m}^{(k)}(x) \rho_{\Xi_n'}^{(k)}\left( y(z,x) \right) \dd x.
\end{align*}
Using the formulas for the correlation functions and the permutational definition of the determinant, we can put the integrand in the following form
\begin{align*}
&\frac{\1_{V_\eta}(x)}{(2\pi)^k}\cdot \det\left[ \frac{2\pi}{m}s_m(x_s - x_t)\right]_{s,t=1}^k \cdot \det\left[ \frac{2\pi}{n}s_n\left( y_s - y_t \right) \right]_{s,t=1}^k \\
&= \frac{\1_{V_\eta}(x)}{(2\pi)^k}\Bigg(1 + \sum_{\sigma\neq\textrm{id} \ \textrm{or} \ \tau\neq\textrm{id}}\sgn \sigma \sgn \tau \prod_{i=1}^k \frac{2\pi}{m}s_m(x_i - x_{\sigma(i)}) \cdot  \frac{2\pi}{n}s_n(y_i - y_{\tau(i)})\Bigg),
\end{align*}
where the second summation runs through permutations $\sigma$ and $\tau$ of k indices. The point is that each term in this sum tends to zero with $m$ and $n$ going to infinity as we have $\frac{2\pi}{m}s_m(x_i - x_{\sigma(i)}) \xrightarrow[m\to \infty]{\textrm{a.e.}} 0$ for $i$ such that $i \neq \sigma(i)$, and $\frac{2\pi}{n}s_n(y_i - y_{\tau(i)}) \xrightarrow[n\to \infty]{\textrm{a.e.}} 0$ if $i \neq \tau(i)$ (see \eqref{eq:simpleobserv} and mind the fact that actually $y$ depends on $m$ and $n$). Recall also that $\1_{W_\eta} \to 1$ and $\sum_\eta \1_{V_\eta} \to \1_{[0,2\pi)^k}$. Moreover, \eqref{eq:easyboundS_n} yields that the whole sum is bounded by $(k!)^2/(2\pi)^k$. Therefore by Lebsegue's dominated convergence theorem we conclude that
\[ \rho_{\Theta_{m,n}}^{(k)}(z) \xrightarrow[m,n \to \infty]{} \int \1_{[0,2\pi)^k}(x)\frac{1}{(2\pi)^k} \dd x = 1, \]
which finishes the proof. \qedsymbol

\section{Concluding remarks}

At the very end we shall discuss the tensor product of more than two matrices. We only briefly sketch what can be easily inferred looking at the proof of the main result.

Let $\Xi_l$, $\Xi_m'$, $\Xi_n''$ be the point processes of eigenphases of independent $l \times l$, $m \times m$,  and $n \times n$ random unitary matrices respectively. Proceeding along the same lines as in the proof of Theorem (\ref{thm(a)}), we conclude that the point process $\frac{2\pi}{lmn}(\Xi_l \otimes \Xi_m' \otimes \Xi_n'' - \pi)$ locally behaves as the Poisson point process on $\R$ when $l$ is fixed but $m$ and $n$ tend to infinity. Indeed, the asymptotics of the $k$-th correlation function $\rho^{(k)}(z)$ of that process is governed by the integrals
\[ \int_{[0,2\pi]^{p+k}\cap V_\eta} \left(\frac{2\pi}{lmn}\right)^{k}\rho^{(p)}_{\Xi_l}(x)\rho^{(k)}_{\Xi_m'}(y)\rho^{(k)}_{\Xi_n''}(w(x,y,z)) \dd x \dd y, \]
which we then sum suitably. Expanding the determinantal correlation functions of $\Xi_m'$ and $\Xi_n''$ (see the proof of Theorem (\ref{thm(b)})) we find that the limit of $\rho^{(k)}(z)$ equals $\sum_{p, \pi} \frac{1}{l^k}\frac{l!}{(l-p)!} = 1$, where the last identity is due to the well-known combinatorial fact that $\sum_{p=1}^k \sharp \mathfrak{S}(k, p) x(x-1)\cdot\ldots(x-p+1) = x^k$. The same line of reasoning applies also when in addition $l \to \infty$. Then the asymptotics depends only on the integral 
\[ \int_{[0,2\pi]^{2k}} \left(\frac{2\pi}{lmn}\right)^{k} \rho^{(k)}_{\Xi_l}(x)\rho^{(k)}_{\Xi_m'}(y)\rho^{(k)}_{\Xi_n''}(w(x,y,z)) \dd x \dd y. \]
Again, we carry on as in the proof of Theorem (\ref{thm(b)}).

Let $A^{(i)}_{n_i}$, $i = 1, ,2, \ldots$, be independent $n_i \times n_i$ random unitary matrices. The other cases of tensor products $\bigotimes_{i=1}^M A^{(i)}_{n_i}$, when for instance all but one of $n_i$'s are fixed, seem to be more delicate and we do not wish to go into detail here. Moreover, it looks challenging to consider the tensor products when the number of terms $M$ tends to infinity and $(n_i)_{i=1}^\infty$ is fixed. The simplest case of $n_i = 2$, $i \geq 1$ has been addressed in \cite{T}.

\section*{Acknowledgements}

Thanks to Prof. Neil O'Connell for his great help.

\end{document}